\documentclass{amsart}

\usepackage{amsmath}
\usepackage{amssymb}
\usepackage{a4wide}
\usepackage{epsfig}

\newcommand{\R}{\ensuremath{\mathbb{R}}}

\renewcommand{\H}{\ensuremath{\mathbb{H}}}
\newcommand{\Sb}{\ensuremath{\mathbb{S}}}

\newcommand{\Sc}{\ensuremath{{\mathcal S}}}

\renewcommand{\rho}{\varrho}

\renewcommand{\phi}{\varphi}

\DeclareMathOperator{\ident}{id}

\newtheorem{thm}{Theorem}[section]
\newtheorem{lemma}[thm]{Lemma}

\begin{document}

\title{Counting perfect colourings of plane regular tilings}

\author{Dirk Frettl\"oh}
\address{Fakult\"at f\"ur Mathematik, Universit\"at Bielefeld, 
Postfach 100131, 33501 Bielefeld,  Germany}
\email{dirk.frettloeh@math.uni-bielefeld.de}
\urladdr{http://www.math.uni-bielefeld.de/baake/frettloe}

\begin{abstract}
A first step in investigating colour symmetries
of periodic and nonperiodic patterns is determining the number of
colours which allow perfect colourings of the pattern under
consideration. A perfect colouring is one where each symmetry of the
uncoloured pattern induces a global permutation of the colours. 
Two cases are distinguished: Either perfect colourings with respect to
all symmetries, or with respect to orientation preserving symmetries
only (no reflections). For the important class of colourings of
regular tilings (and some Laves tilings) of the Euclidean or hyperbolic
plane, this mainly combinatorial question is addressed here using
group theoretical methods.
\end{abstract}

\maketitle

\section{ \bf Introduction}
The study of colour symmetries of crystallographic patterns in Euclidean
space is a classical topic, see for instance \cite{gs}, \cite{schw}
and references therein. The discovery of quasicrystals inspired
interest in colour symmetries of further (quasi-)crystallographic
patterns, see \cite{baa}, \cite{lip}, \cite{mp}. A very recent  
development is the study of colour symmetries of regular tilings in
the hyperbolic plane \cite{mlp}. For several reasons, the focus
of studying colour symmetries frequently lies on perfect colourings
(for a precise definition, see below). For instance, a perfect
colouring yields a simple algebraic relation between the symmetry
group $G$ of the uncoloured pattern and the symmetry group $K$ of the
coloured pattern: The factor group $G/K$ is a permutation group. 
Thus, basic questions to address in the study of colour symmetries are
determining the possible numbers of colours for which a perfect
colouring exists, counting the different perfect colourings with this 
number of colours, and study algebraic properties of the colour
groups. This article answers the first two questions for perfect
colourings of regular (and certain Laves tilings), be it spherical
tilings, or tilings in the Euclidean plane, or in the hyperbolic
plane. The case of Euclidean tilings is well known for long, see
\cite{schw} for a good survey. The case of regular tilings of the
sphere is somehow trivial. Therefore, the focus of this article lies
on counting perfect colourings of hyperbolic regular (and certain
hyperbolic Laves) tilings. The main result is a method to count these
perfect colourings based on counting certain subgroups of Coxeter
groups. This is outlined in Section \ref{zwei} with respect to the
full symmetry group of regular tilings, and in Section \ref{drei} with
respect to the rotation group of regular tilings. As an application,
some concrete results are listed in Table \ref{tabelle1} and
\ref{tabelle2}. 

\section{ \bf Colourings, regular tilings and Laves tilings} 
Let $\mathbb{X}$ be either the Euclidean plane
$\R^2$, or the hyperbolic plane $\H^2$, or the sphere $\Sb^2$. An
isometry $f : \mathbb{X} \to \mathbb{X}$ is called a {\em symmetry}
of some set $X \subset \mathbb{X}$, if $f(X)=X$. Here, $X$ can be
nearly anything living in $\mathbb{X}$, but it is fine to think of $X$
as some point set, or some tiling. The set of all symmetries of $X$
form a group, the {\em symmetry group} of $X$. A symmetry is called 
{\em orientation-preserving}, if it is a translation, or a rotation,
or a product of those. The subgroup of the symmetry group of $X$ 
containing all orientation-preserving symmetries is called the 
{\em rotation group} of $X$. In the remainder of this section, the
symmetry group of some set $X$ is always denoted by $G$, as well as
the rotation group of $X$. This should not lead to confusion, 
everything works for both cases.

In the sequel, the symmetric group of order $n$ (all permutations on
$n$ numbers) is denoted by $\Sc_n$. A {\em k-colouring} $(X,c)$ of $X$
is a surjective map $c: X \to \{1, \ldots, k\}$. We are mainly
interested in colourings, where each $h \in G$ acts as a global
permutation of colours. That is, we consider those $h \in G$ where 
\[ \forall x,y \in X: \; c(x)=c(y) \Rightarrow c(h(x))=c(h(y)). \]
Following \cite{gs}, we call a pair $(h, \pi)$ (where $h \in G,
\pi \in \Sc_k$) a {\em colour symmetry} of $(X,c)$, if $\pi (c(x)) =
c(h(x))$ for all $x \in X$. The group of all colour symmetries  
\[ C(X) = \{ (h, \pi) \; | \; h \in G,  \pi \in \Sc_k : c(h(x)) =
\pi (c(x))\; \mbox{for all} \; x \in X \} \]
is called the {\em colour symmetry group} of $X$.  
For convenience, let $H := \{ h \; | \; (h,\pi) \in C(X) \}$. 
Since each $h \in H$ determines a unique permutation $\pi_h$, $H$ is
isomorphic to $C(X)$. The map $h \mapsto \pi_h$ is a group
homomorphism, which can be seen as follows: Let $i \in \{1,
\ldots,k\}$. Then, for all $x \in X$ with $c(x)=i$ holds 
\[ \pi_{gh}(i) = \pi_{gh}(c(x)) = c(gh(x)) = c(g(h(x))) = 
\pi_g(c(h(x)) = \pi_g \pi_h (c(x)) = \pi_g \pi_h (i), \]
hence $\pi_{gh}=\pi_g \pi_h$.   
A colouring is called {\em perfect}, if $H=G$, or
equivalently: if every symmetry of $X$ is a colour symmetry of
$(X,c)$. 

\begin{figure}[hbt]
\begin{center}
\epsfig{file=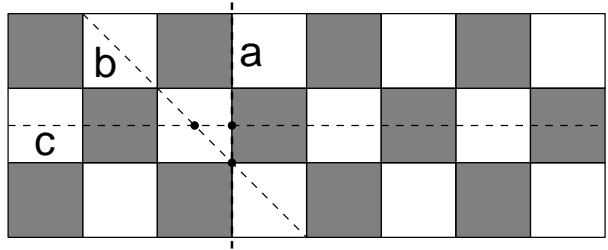, width=80mm}
\end{center}
\caption{Part of an infinite checkerboard. A reflection in the line
  $a$ interchanges the colours globally. Reflections in the lines $b$
  or $c$ fix the colours. \label{fig:schach}} 
\end{figure}

{\em Example:} Consider an infinite checkerboard (Figure
\ref{fig:schach}). It can be regarded as a 2-colouring of the
canonical square tiling of the plane. The symmetry group $G$ of the
uncoloured tiling is the Coxeter group generated by reflections in
the lines $a,b,c$. A reflection $R_a$ in $a$ interchanges the colours
of the checkerboard: 
All black squares are mapped to white squares and vice versa. Thus,
$R_a \in H$, $\pi_{R_a} = (1,2)$. The reflections $R_b$ and $R_c$ both
fix the colours, therefore $R_b, R_c \in H$, $\pi_{R_b}=\pi_{R_c} =
\ident$. All three generators of $G$ are also elements of $H$, thus
$G=H$, and the colouring is perfect. 

A tiling of $\mathbb{X}$ is a collection of compact sets
$T_1, T_2, \ldots$ which cover $\mathbb{X}$ (that is, the union of all $T_i$
is $\mathbb{X}$) and which do not overlap (that is, the intersection of the
interiors of any two sets $T_i \ne T_j$ are empty). The sets $T_i$ are
called tiles. A tiling of $\mathbb{X}$ is called {\em regular}, if
all tiles are regular polygons, and all vertex figures are regular
polygons. (The vertex figure of a vertex $x$ in a tiling is the convex
hull of all midpoints of the edges emanating from $x$.) Following
\cite{gs}, a tiling of $\R^2$ or $\H^2$ is called {\em Laves tiling}, 
if all tiles are congruent, and all vertex figures are regular
polygons. For a thorough survey on regular tilings and Laves tilings see
\cite{cox} and \cite{gs}. 

Any ordered pair of integers $(p,q)$, where $p,q \ge 3$, defines a
regular tiling by regular $p$-gons, where $q$ tiles meet at 
each vertex. Following \cite{gs}, such a tiling is denoted by
$(p^q)$.  If $\frac{1}{p} + \frac{1}{q} > \frac{1}{2}$, then
$(p^q)$ is a  regular tiling of the sphere (which can be regarded as
a regular polytope in $\R^3$). If $\frac{1}{p} + \frac{1}{q} =
\frac{1}{2}$, then $(p^q)$ is a  regular tiling of the Euclidean
plane $\R^2$.  If $\frac{1}{p} + \frac{1}{q} < \frac{1}{2}$,  
then $(p^q)$ is a regular tiling of the hyperbolic plane $\H^2$. 
In $\R^2$ there are only three different regular tilings and eleven
different Laves tilings (including the three regular ones, see
\cite{gs}). In $\Sb^2$ there are five regular tilings (corresponding
to the five platonic solids) and eight Laves tilings (including the
former five, see \cite{cox}). In contrast, there are infinitely many
regular tilings in $\H^2$, as well as Laves tilings. Here we consider
the following Laves tilings only: Each regular tiling $(p^q)$ gives
rise to a Laves tiling $[p.q.p.q]$ (if $p<q$, otherwise we write
$[q.p.q.p]$, in order to make the notation unique) by joining the
centres of each $p$-gons with each of its vertices by line
segments. Consider these line segment as the edges of a new
tiling. Then this tiling will be a tiling by congruent quadrilaterals,
with either $p$ or $q$ tiles at each vertex. (Note that for the
regular square tiling in $\R^1$, this constructions yields again a
regular square tiling $[4.4.4.4] = (4^4)$; whereas the regular hexagon
tiling $(6^3)$ yields a Laves tiling $[3.6.3.6]$ by $60^{\circ}$
rhombi.)  

The symmetry group of $(p^q)$ is a Coxeter group $G_{p,q}$:
\begin{equation} \label{gpq}
G_{p,q} = \langle a,b,c \, | \, a^2=b^2=c^2= (ab)^q 
= (ac)^2 = (bc)^p = 1 \rangle. 
\end{equation} 
$G_{p,q}$ is also the symmetry group of both $(q^p)$ and 
$[p.q.p.q]$. The latter fact is the reason, why we can easily include
these certain Laves tilings in our study, where other Laves tilings
require further effort.
A fundamental  domain of $G_{p,q}$ is an orthogonal triangle $F_{p,q}$ 
with angles $\frac{\pi}{2}, \frac{\pi}{p}, \frac{\pi}{q}$. Then, $a$
($b,c$, resp.) denotes the reflection in the line spanned by the edge
of $F$ which is opposite to the angle $\frac{\pi}{p}$ ($\frac{\pi}{2},
\frac{\pi}{q}$, resp.).  

{\bf Left coset colouring:} (compare \cite{mlp})
The orbit of $F:=F_{p,q}$ under $G_{p,q}$ is a tiling $T_{p,q}$ of
$\mathbb{X}$ by  triangles $f F$, where $f$ runs through $G_{p,q}$. 
In particular, this gives a bijection between elements of $G_{p,q}$
and elements of $T_{p,q}$. Therefore, any subgroup $S \subseteq
G_{p,q}$ of index $k$ induces a $k$-colouring of $T_{p,q}$ in the
following way: Set $c(f F)=1$ for each $f \in S$, and set $c(f F)= i$,
if $f$ is element of the $i$-th left coset $S_i$ of $S=S_1$ ($1 \le i
\le k$). In the next section, we deduce colourings of regular tilings
$(p^q)$ from these colourings of $T_{p,q}$. 

This construction is taken from \cite{mlp}, where it is given in a
more general setup. Here it reduces to the above, since the stabiliser
of $F$ in $G_{p,q}$ is trivial. However, the following lemma is true
in a more general context than that of this article.

\begin{lemma} \label{perfect}
All colourings of a tiling obtained by the left coset colouring method
are perfect.  
\end{lemma}
{\em Proof:}
Let $S$ be the subgroup of $G_{p,q}$ inducing the colouring of
$T_{p,q}$. Each tile in $T_{p,q}$ is of the form $fF$, where $f \in
G_{p,q}$.  We need to show that for each symmetry $g \in G_{p,q}$
holds: $c(fF) = c(f'F)$ implies $c(gfF)=c(gf'F)$. By the left coset
construction, $c(fF)=c(f'F)$ is equivalent to $f$ and $f'$ being in
the same left coset. That is, $f^{-1} f' \in S$. Consequently,
$(gf)^{-1}(gf') = f^{-1}g^{-1}gf' = f^{-1} f' \in S$, and the claim
follows.  \hfill $\square$

\section{ \bf Counting colourings with respect to the symmetry group}
\label{zwei}  
In this section we present a method to count colour symmetry groups
for any regular tiling or Laves tiling in $\R^2$ or $\H^2$ or $\Sb^2$
with respect to the entire symmetry group of the tiling.
We should emphasise that different colourings
may possess the same colour symmetry group (see \cite[\S 8]{schw}). 
Colourings with the same colour symmetry group are strongly 
related, but not necessarily congruent. So what we count here are not
different colourings, but different colour symmetry groups, where
different means non-conjugate. This is the usual approach, compare
\cite{bg}, \cite{schw}. In general, the question how two different
colourings with the same colour symmetry group are related can be
quite hard to answer. 

Note that $T_{p,q}$ can be obtained from $(p^q)$ by dividing each
$p$-gon into $2p$ triangles congruent to $F$. In the following, we are 
interested in colourings of $(p^q)$ (respectively $[p.q.p.q]$,
respectively $(q^p)$) rather than colourings of $T_{p,q}$. To obtain a
colouring (of $(p^q)$, say) we just need to consider those colourings
of $T_{p,q}$, where all $2p$ triangles forming a $p$-gon possess the
same colour.  

\begin{lemma}
In a perfect colouring of a regular tiling $(p^q)$ or a Laves tiling
$[p.q.p.q]$, all colours occur with the same frequency.
\end{lemma}
{\em Proof:}
The symmetry group $G_{p,q}$ acts transitively on the tiles. Thus, for all
pairs of colours $(i,j)$, there is $f_{ij} \in G_{p,q}$ mapping some tile
of colour $i$ to some tile of colour $j$. Since the colouring is perfect,
$f_{ij}$ maps the entire colour class $c^{-1}(i)$ to the colour class
$c^{-1}(j)$. Since $f_{ij}$ is an isometry, $c^{-1}(i)$ and $c^{-1}(j)$ 
have the same frequency. \hfill $\square$

The last result is not true in general for perfect colourings of other
tilings.  

\begin{lemma} \label{subgr}
Let $S$ be a subgroup of $G_{p,q}$ of index $k$. 
$S$ induces a perfect $k$-colouring of $(q^p)$, if and only if 
$a,b \in S$. 
$S$ induces a perfect $k$-colouring of $(p^q)$, if and only if 
$b,c \in S$. 
$S$ induces a perfect $k$-colouring of $[p.q.p.q]$, if and only if 
$a,c \in S$.
\end{lemma}
{\em Proof:}
We provide the proof for the case $(q^p)$. The other cases are
completely analogous. 

One direction ('if') is clear from the construction and from Lemma
\ref{perfect} above. 

For the other direction, let $((q^p),c)$ be a perfect $k$-colouring. 
Let $K_i = \{ f \in G_{p,q} \, | \, c(fF)= i \}$. Without loss of
generality, let $\ident \in K_1$. Then 
\[ \pi(K_1) = \{ \pi_f \, | \, f \in K_1 \} = \{ \pi_f \, | \,
\pi_f(1)=1 \} \] 
is a subgroup of the symmetric group $\Sc_k$. Since $f \mapsto \pi_f$
is a group homomorphism, $K_1$ is a subgroup of $G_{p,q}$. By the
last lemma, all colour classes have the same frequency, thus the index
of $K_1$ in $G_{p,q}$ is $k$. Since the entire $q$-gon $F \cup aF \cup
abF \cup abaF\ldots$ has colour 1, one has $a,ab,aba, \ldots \in K_1$.
Finally, let $i$ and $fF \in K_i$ be fixed. For each
$gF \in K_i$ holds $c(gF)=i$, hence $\pi_{g}(1)=i$, and
$(\pi_g)^{-1} (i) = \pi_{g^{-1}} (i)= 1$. It follows
\[ c(g^{-1}fF) = \pi_{g^{-1}} c(fF) = \pi_{g^{-1}} \pi_f (c(F)) =  
\pi_{g^{-1}} \pi_f (1) = \pi_{g^{-1}} (i) = 1. \]
Therefore, $g^{-1}fF \in K_1$ for all $g \in K_i$, thus $K_i = f K_1$
is a coset of $K_1$. \hfill $\square$

By the last result, we obtain all perfect $k$-colourings of
$(q^p)$ (respectively $[p.q.p.q]$, respectively $(p^q)$) by listing
the subgroups of $G_{p^q}$ which contain $a,b$ (respectively $a,c$,
respectively $b,c$). This can be carried out by GAP \cite{gap}, for
instance. See Table \ref{tabelle1} for some examples. Subgroups are
identified if they are conjugate in $G_{p,q}$.

\begin{table}[h]
\begin{tabular}{|c|l|}
\hline 
$(7^3)$ &  $1, 8, 15, 22, 24, 30, 36^{2}, 44, 50^{5}, \ldots$ \\
\hline     
$[3.7.3.7]$ & $1, 9, 15, 21, 28^{2}, 30, 35^{2}, 36, 37, 42^{5}, 49^{8},50^{3}, \ldots$\\ 
\hline     
$(3^7)$ & $1, 22, 28^{5}, 37, 42^{4}, 44, 49^{7}, 50^{3}, \ldots$\\
\hline     
\hline     
$(8^3)$ &  $1, 3, 6, 12, 17, 21^{4}, 24, 25^{5}, 27^{3}, 29^{4},31^{4}, 33^{6}, 37^{6}, 39^{8}, \ldots$ \\ 
\hline     
$[3.8.3.8]$ & $1, 3, 6, 12^{3}, 17, 18, 21^{4}, 24^{15}, 25^{5},27^{3}, 28^{4}, 29^{4}, 30^{7}, \ldots$ \\
\hline     
$(3^8)$ & $1, 2, 4, 8, 10^{2}, 12, 14, 16^{2}, 18, 20^{4}, 24^{3},25^{5}, 26, 28^{12}, 29, 30^{2}, \ldots$\\ 
\hline     
\hline     
$(5^4)$  & $1, 2, 6, 11, 12, 16^{2}, 21^{3}, 22^{5}, 24, 26^{9}, 28,\ldots$ \\ 
\hline     
$[4.5.4.5]$ & $1, 5^{2}, 10^{5}, 11, 15^{7}, 20^{22}, 21^{3}, 22^{4},25^{27}, 26^{4}, 27^{3}, 30^{63}, \ldots$\\ 
\hline     
$(4^5)$ & $1, 5^{2}, 10^{4}, 11, 15^{7}, 16, 20^{9}, 21^{3}, 22,25^{27}, 26, 27^{3}, 30^{38}, \ldots$\\ 
\hline 

\end{tabular}
\caption{The list of the first possible values $k$ for perfect
  $k$-colourings of some hyperbolic regular and Laves tilings.
  An upper index $n \ge 2$ means that there are $n$ non-conjugate
  colour symmetry groups for this value. No upper index means there is
  exactly one colour symmetry group for this value, up to
  conjugacy. \label{tabelle1}} 
\end{table}

The procedure is quite general. In principle, all colour symmetries 
of all plane regular tilings can be obtained in this
manner, be it hyperbolic, Euclidean or spherical tilings. Regular
spherical tilings in $\Sb^2$ can be regarded as regular polytopes in
$\R^3$. For instance, the cube $(4^3)$ allows perfect colourings
with 1, 3 and 6 colours, each one unique, up to permutation of
colours. The icosahedron $(3^5)$ allows unique perfect colourings with
1, 10 and 20 colours. However, this method has a drawback:
Practically, computing with finitely presented groups is not
efficient. Therefore, there is for each $(p^q)$ a limit on the
number of colours to be considered. Table \ref{tabelle1} lists the
first values of the number of perfect colourings of some regular and
some related Laves tilings in $\H^2$.

\section{ \bf Counting colourings with respect to the rotation group}
\label{drei}  

In the last section colourings were considered which are perfect with
respect to the entire symmetry group $G_{p,q}$. A slightly weaker, but
still important concept is to consider colourings which are perfect
with respect to the rotation group. These include of course all
colourings of the last section, but also many pairs of enantiomorphic
colourings, where one is obtained from the other by a reflection.  
Let us denote the rotation group of the regular tiling $(p^q)$ by
$\bar{G}_{p,q}$. It is also the symmetry group of both $(q^p)$ and
$[p.q.p.q]$. It is a subgroup of $G_{p,q}$ of index 2. $\bar{G}_{p,q}$
is spanned in $G_{p,q}$ by 
$ab, ac$, compare \eqref{gpq}. (With $ab$ and $ac$ it contains also
$bc$. It does not matter which two generators out of $\{ab, ac, bc\}$
are chosen.) Now, analogously to Lemma \ref{subgr}, certain subgroups
of $\bar{G}_{p,q}$ induce colourings of $(p^q)$ or $[p.q.p.q]$ which
are perfect with respect to $\bar{G}_{p,q}$.    

\begin{lemma}
Let $S$ be a subgroup of $\bar{G}_{p,q}$ of index $k$. 
$S$ induces a perfect $k$-colouring of $(q^p)$ if and only if 
$ab \in S$. 
$S$ induces a perfect $k$-colouring of $[p.q.p.q]$ if and only if 
$ac \in S$. 
$S$ induces a perfect $k$-colouring of $(p^q)$ if and only if 
$bc \in S$.
\end{lemma}
Note that 'perfect' is to be read here as 'perfect with respect to the
rotation group'. The proof is completely analogous to the one of Lemma
\ref{subgr}, using the fact, that $F \cup aF$ and $F \cup bF$ are both
fundamental domains of $\bar{G}_{p,q}$.  

Counting the perfect colourings is slightly more subtle as in the last
section. Again subgroups are identified if they are conjugate.
But two subgroups have to be
identified not only if they are conjugate in $\bar{G}_{p,q}$, but also
if they are conjugate in $G_{p,q}$. This requires some programming in
GAP, and the computations become even harder. (The program code can be
requested from the author.) For instance, there are
119 subgroups of index 20 in $\bar{G}_{5,4}$ inducing perfect
colourings of $[4.5.4.5]$. These have to be checked pairwise for being
conjugate in $G_{5,4}$. Without further refinements, this takes weeks
using GAP on a usual PC, which is the reason why this particular value
has not been included in Table \ref{tabelle2}. The
groups in Table \ref{tabelle1} are contained in Table \ref{tabelle2}.

\begin{table}[h]
\begin{tabular}{|c|l|}
\hline 
$(7^3)$ & $1, 8, 9, 15^{2}, 22^{7}, 24,  \ldots$ \\
\hline     
$[3.7.3.7]$ & $1,  7, 9, 14^{6}, 15^{2}, 21^{5}, 22^{7}, \ldots$\\ 
\hline     
$(3^7)$ & $1,  7, 8, 14^{6}, 21^{2}, 22^{7}, \ldots$\\
\hline     
\hline     
$(8^3)$ & $1,  3, 6, 9, 10, 12, 13^{2}, 15, 17^{5}, 18^{5}, 19^{5}, \ldots$ \\  
\hline     
$[3.8.3.8]$ & $1,  3, 4, 6, 8^{3}, 9, 10^{2}, 12^{4}, 13^{2}, 15, 16^{9}, 17^{5}, 18^{12}, \ldots$ \\
\hline     
$(3^8)$ & $1,  2, 4, 8^{4}, 10^{3}, 12, 13^{2}, 14^{2}, 16^{12}, 17^{5}, 18, 19^{5}, \ldots$\\ 
\hline     
\hline     
$(5^4)$  & $1, 2, 6^{2}, 11^{3}, 12^{6}, 16^{12}, 17^{4}, \ldots$ \\ 
\hline     
$[4.5.4.5]$ & $1, 5^{2}, 6, 10^{9}, 11^{3}, 12^{4}, 15^{15}, 16^{10},17^{4},  \ldots$\\  
\hline     
$(4^5)$ & $1, 5^{2}, 6, 10^{6}, 11^{3}, 15^{15}, 16^{2}, 17^{4},\ldots$\\  
\hline  

\end{tabular}
\caption{The list of the first possible values $k$ of perfect
  $k$-colourings of some hyperbolic regular and Laves tilings, 
  using the same scheme as Table \ref{tabelle1}. The groups in
  Table \ref{tabelle1} are contained in this table. \label{tabelle2}} 
\end{table}

\section{ \bf Conclusion}
The left coset colouring method of \cite{mlp} to obtain perfect
colourings of certain tilings was used to obtain perfect colourings of
regular tilings $(p^q)$ (or Laves tilings $[p.q.p.q]$). The use of
Coxeter groups $G_{p,q}$ allowed for a bijection between perfect
colourings of these tilings and subgroups of $G_{p,q}$. This allowed
us to deduce the numbers $k$ for which a $k$-colouring exists, and
their multiplicities (up to conjugacy of the colour symmetry
groups). As an application, we provide a list of these numbers $k$
together with their multiplicities for some hyperbolic 
regular and Laves tilings. 

\section{ \bf Acknowledgements}
It is a pleasure to thank the second referee for very
valuable remarks. This work was supported by the German
Research Council (DFG) within the CRC 701.

\end{document}